\documentclass[twoside,leqno,10pt, A4]{amsart}
\usepackage{amsfonts}
\usepackage{amsmath}
\usepackage{amscd}
\usepackage{amssymb}
\usepackage{amsthm}
\usepackage{amsrefs}
\usepackage{latexsym}
\usepackage{mathrsfs}
\usepackage{bbm}
\usepackage{amscd}
\usepackage{amssymb}
\usepackage{amsthm}
\usepackage{amsrefs}
\usepackage{latexsym}
\usepackage{mathrsfs}
\usepackage{bbm}
\usepackage{enumerate}
\usepackage{graphicx}
\usepackage{color}
\begin{document}

\newtheorem{theorem}[subsection]{Theorem}
\newtheorem{proposition}[subsection]{Proposition}
\newtheorem{lemma}[subsection]{Lemma}
\newtheorem{corollary}[subsection]{Corollary}
\newtheorem{conjecture}[subsection]{Conjecture}
\newtheorem{prop}[subsection]{Proposition}
\newtheorem{defin}[subsection]{Definition}

\numberwithin{equation}{section}
\newcommand{\mr}{\ensuremath{\mathbb R}}
\newcommand{\mc}{\ensuremath{\mathbb C}}
\newcommand{\dif}{\mathrm{d}}
\newcommand{\intz}{\mathbb{Z}}
\newcommand{\ratq}{\mathbb{Q}}
\newcommand{\natn}{\mathbb{N}}
\newcommand{\comc}{\mathbb{C}}
\newcommand{\rear}{\mathbb{R}}
\newcommand{\prip}{\mathbb{P}}
\newcommand{\uph}{\mathbb{H}}
\newcommand{\fief}{\mathbb{F}}
\newcommand{\majorarc}{\mathfrak{M}}
\newcommand{\minorarc}{\mathfrak{m}}
\newcommand{\sings}{\mathfrak{S}}
\newcommand{\fA}{\ensuremath{\mathfrak A}}
\newcommand{\mn}{\ensuremath{\mathbb N}}
\newcommand{\mq}{\ensuremath{\mathbb Q}}
\newcommand{\half}{\tfrac{1}{2}}
\newcommand{\f}{f\times \chi}
\newcommand{\summ}{\mathop{{\sum}^{\star}}}
\newcommand{\chiq}{\chi \bmod q}
\newcommand{\chidb}{\chi \bmod db}
\newcommand{\chid}{\chi \bmod d}
\newcommand{\sym}{\text{sym}^2}
\newcommand{\hhalf}{\tfrac{1}{2}}
\newcommand{\sumstar}{\sideset{}{^*}\sum}
\newcommand{\sumprime}{\sideset{}{'}\sum}
\newcommand{\sumprimeprime}{\sideset{}{''}\sum}
\newcommand{\sumflat}{\sideset{}{^\flat}\sum}
\newcommand{\shortmod}{\ensuremath{\negthickspace \negthickspace \negthickspace \pmod}}
\newcommand{\V}{V\left(\frac{nm}{q^2}\right)}
\newcommand{\sumi}{\mathop{{\sum}^{\dagger}}}
\newcommand{\mz}{\ensuremath{\mathbb Z}}
\newcommand{\leg}[2]{\left(\frac{#1}{#2}\right)}
\newcommand{\muK}{\mu_{\omega}}
\newcommand{\thalf}{\tfrac12}
\newcommand{\lp}{\left(}
\newcommand{\rp}{\right)}
\newcommand{\Lam}{\Lambda_{[i]}}
\newcommand{\lam}{\lambda}
\newcommand{\af}{\mathfrak{a}}
\newcommand{\sw}{S_{[i]}(X,Y;\Phi,\Psi)}
\newcommand{\lz}{\left(}
\newcommand{\pz}{\right)}
\newcommand{\bfrac}[2]{\lz\frac{#1}{#2}\pz}
\newcommand{\odd}{\mathrm{\ primary}}
\newcommand{\even}{\text{ even}}
\newcommand{\res}{\mathrm{Res}}
\newcommand{\sumn}{\sumstar_{(c,1+i)=1}  w\left( \frac {N(c)}X \right)}
\newcommand{\lab}{\left|}
\newcommand{\rab}{\right|}
\newcommand{\Go}{\Gamma_{o}}
\newcommand{\Ge}{\Gamma_{e}}
\newcommand{\M}{\widehat}

\theoremstyle{plain}
\newtheorem{conj}{Conjecture}
\newtheorem{remark}[subsection]{Remark}

\makeatletter
\def\widebreve{\mathpalette\wide@breve}
\def\wide@breve#1#2{\sbox\z@{$#1#2$}%
     \mathop{\vbox{\m@th\ialign{##\crcr
\kern0.08em\brevefill#1{0.8\wd\z@}\crcr\noalign{\nointerlineskip}%
                    $\hss#1#2\hss$\crcr}}}\limits}
\def\brevefill#1#2{$\m@th\sbox\tw@{$#1($}%
  \hss\resizebox{#2}{\wd\tw@}{\rotatebox[origin=c]{90}{\upshape(}}\hss$}
\makeatletter

\title[Ratios conjecture of quadratic Hecke $L$-functions of prime-related moduli]{Ratios conjecture of quadratic Hecke $L$-functions of prime-related moduli}

%%\date{\today}
\author[P. Gao]{Peng Gao}
\address{School of Mathematical Sciences, Beihang University, Beijing 100191, China}
\email{penggao@buaa.edu.cn}

\author[L. Zhao]{Liangyi Zhao}
\address{School of Mathematics and Statistics, University of New South Wales, Sydney NSW 2052, Australia}
\email{l.zhao@unsw.edu.au}

\begin{abstract}
Using the method of multiple Dirichlet series, we develop $L$-functions ratios conjecture with one shift in both the numerator
and denominator in certain ranges for quadratic families of Dirichlet and Hecke $L$-functions of prime-related moduli of imaginary quadratic number fields of class number one under the generalized Riemann hypothesis. As corollaries, we evaluate asymptotically the first moment of central values as well as the one-level density of the families of $L$-functions under consideration. 
\end{abstract}

\maketitle



\section{Introduction}
\label{sec 1}

Originating in the work of D. W. Farmer \cite{Farmer93} concerning the shifted moments of the Riemann zeta function and formulated for general $L$-functions in the work of J. B. Conrey, D. W. Farmer and M. R. Zirnbauer \cite[Section 5]{CFZ}, the $L$-functions ratios conjecture predicts the asymptotic behaviors of the sum of ratios of products of shifted $L$-functions. This conjecture has many important number theoretic applications, most notably in the density conjecture of N. Katz and P. Sarnak \cites{KS1, K&S} on the distribution of zeros near the central point of a family of $L$-functions. \newline

 The ratios conjecture for quadratic families of $L$-functions was first developed by H. M. Bui, A. Florea and J. P. Keating  \cite{BFK21} over function fields and later by M. \v Cech \cite{Cech1} over $\mq$. The authors further studied the ratios conjecture for various families of $L$-functions in \cites{G&Zhao14, G&Zhao15, G&Zhao16, G&Zhao17}.  \newline 
 
  In this paper, we are interested in investigating the ratios conjecture for families of quadratic Hecke $L$-functions of prime moduli. These $L$-functions differ from those studied in \cite{BFK21} and \cite{Cech1} as their conductors form only sparse sets (sets of density zero) in the relevant ring of integers. \newline
  
  To state our results, suppose that $K$ is an imaginary quadratic number field of class number one or the field of rational numbers $\mq$.  Hence $K=\mq$ or $K=\mq(\sqrt{d})$ (see \cite[(22.77)]{iwakow}) with $d \in \mathcal{S}$, where
\begin{align*}
%%\label{dvalue}
\mathcal{S} = \{-1, -2, -3, -7,-11,-19,-43,-67,-163 \}.
\end{align*}
Let $\mathcal{O}_K, U_K$ and $D_K$ denote the ring of integers, the group of units and the discriminant of $K$, respectively.  We shall focus on these number fields in this paper and note here that for number fields with class number greater than one, ideals in their rings of integers become much harder to describe. \newline

   In what follows, we shall always assume that $K$ is of the above form if it is an imaginary quadratic number field. We then define $c_K=(1+i)^5$ when $d=-1$, $c_K=4\sqrt{-2}$ when $d=-2$ and $c_K=8$ for the other $d$ in $\mathcal{S}$ as well as when $K=\mq$.  Let $L(s, \chi)$ denote as usual the $L$-function attached to any Hecke character $\chi$ and $\zeta_K(s)$ for the Dedekind zeta function of $K$.  Moreover, let $r_K$ denote the residue of $\zeta_K(s)$ at $s = 1$.  We also use the notation $L^{(c)}(s, \chi)$ for the Euler product defining $L(s, \chi)$ but omitting those primes dividing $c$. \newline
   
   Let $\chi^{(c_K\varpi)}$ be the quadratic Hecke character of trivial infinite type defined in Section~\ref{sec2.4}, where we reserve the symbol $\varpi$ for a prime number in $\mathcal{O}_K$, which means that the ideal $(\varpi)$ generated by $\varpi$ is a prime ideal. We are primarily interested in developing the ratios conjecture for quadratic families of Hecke $L$-functions of prime moduli. However, we observe that $\chi^{(\varpi)}$ is not necessarily primitive modulo $\varpi$ for any prime $\varpi$, even in the case when $K=\mq$. Meanwhile, it is shown in \cite[Section 2.1]{G&Zhao4} and \cite[Section 2.2]{G&Zhao2022-4} that $\chi^{(c_K\varpi)}$ is a primitive quadratic character of trivial infinite type for any prime $\varpi$. Thus, we modify our object slightly by considering the ratios conjecture involving with $L(s, \chi^{(c_K\varpi)})$. \newline
   
   We fix a non-negative smooth function $w(t)$ that is compactly supported on ${\mr}^+$, the set of positive real numbers.  We write $N(n)$ for the norm of any $n \in \mathcal{O}_K$ and $\Lambda_K(n)$ for the von Mangoldt function on $\mathcal{O}_K$, given by
\begin{align*}
    \Lambda_K(n)=\begin{cases}
   \log N(\varpi) \qquad & n=\varpi^k, \text{$\varpi$ primary prime}, k \geq 1, \\
     0 \qquad & \text{otherwise}.
    \end{cases}
\end{align*}

  Let $\Gamma(s)$ be the usual gamma function.  We further define
\begin{align} \label{GammaKBKdef}
\begin{split}
  \Gamma_K(s) =
\begin{cases}
 \frac {\Gamma(\frac {1-s}{2})}{\Gamma(\frac {s}{2})}, \quad K =\mq,  \\
  \frac {\Gamma(1-s)}{\Gamma(s)}, \quad \text{otherwise},
\end{cases}
\quad \mbox{and} \quad  B_K =
\begin{cases}
 \frac {|D_K|N(c_K)}{\pi}, \quad K =\mq,  \\
 \frac {|D_K|N(c_K)}{(2\pi)^2}, \quad \text{otherwise}.
\end{cases}
\end{split}
\end{align} 

  Our main result investigates under the generalized Riemann hypothesis (GRH) the ratios conjecture with one shift in both the numerator
and denominator for the families of quadratic Hecke $L$-functions $L(s, \chi^{(c_K\varpi)})$.  Note that one may also obtain similar results for families of quadratic Hecke $L$-functions over general moduli and this has been done in \cite{G&Zhao14} over the Gaussian field $\mq(i)$. In either situation, one may use GRH to better control the error terms. However, one may obtain unconditional results with weaker error terms as well. \newline

  In the sequel, for a non-negative Schwartz function $w$, we use $\widehat{w}$ to denote its Mellin transform so that
\begin{align}
\label{Mellin}
     \widehat{w}(s) =\int\limits^{\infty}_0w(t)t^s\frac {\dif t}{t}.
\end{align}
 Moreover,  $\varepsilon$ always, as is standard, denotes a small positive real number which may not be the same at each occurrence.  
\begin{theorem}
\label{Theorem for all characters}
	With the notation as above and assuming the truth of GRH, let $K=\mq(\sqrt{-d})$ with $d \in \mathcal  S$ or $K=\mq$ and $d_K=[K:\mq]$.  Suppose further that $X$ is a large real number and $w(t)$ a non-negative Schwartz function.  Set
\begin{align}
\label{Nab}
  E(\alpha,\beta)=\max\left\{\tfrac 12, \ 1-\Re(\alpha)-\Re(\beta), \ 1-\tfrac {\Re(\alpha)}2-\tfrac {\Re(\beta)}2,  \ \tfrac {3}{4}-\tfrac {\Re(\alpha)}2, \ 1-\Re(\beta) \right\}.
\end{align}
 Then we have for $-1/4<\Re(\alpha) < 1/2, \ \Re(\beta)>0$ with $\Re(\alpha) +\Re(\beta)>0$,
\begin{align}
\label{Asymptotic for ratios of all characters}
\begin{split}	
\sum_{\substack{\varpi \odd}} & \frac {\Lambda_{K}(\varpi) L(\tfrac 12+\alpha, \chi^{(c_K\varpi)})}{L(\tfrac 12+\beta, \chi^{(c_K\varpi)})}w \bfrac {N(\varpi)}X  \\
=&  X\M w(1)\frac {\zeta^{(2)}_K(1+2\alpha)}{\zeta^{(2)}_K(1+\alpha+\beta)}+X^{1-\alpha}\M w(1-\alpha)B^{-\alpha}_K\Gamma_K(\tfrac 12+\alpha)\frac {\zeta^{(2)}_K(1-2\alpha)}{\zeta^{(2)}_K(1-\alpha+\beta)}  \\ 
& \hspace*{2cm} +O\lz |1+\alpha|^{3-d_K/2+\varepsilon}|1+\beta|^{\varepsilon} X^{E(\alpha,\beta)+\varepsilon}\pz. 
\end{split}
\end{align}
\end{theorem}

The conditions imposed on $\alpha$ and $\beta$ in the statement of Theorem \ref{Theorem for all characters} ensures that  $E(\alpha,\beta)<1$.  Our result above is consistent with the prediction from the ratios conjecture on the left-hand side of \eqref{Asymptotic for ratios of all characters}, which can be derived following the recipe given in \cite[Section 5]{CFZ} except that (see also \cite{CS}) the ratios conjecture asserts that \eqref{Asymptotic for ratios of all characters} holds uniformly for $|\Re(\alpha)|< 1/4$, $(\log X)^{-1} \ll \Re(\beta) < 1/4$ and $\Im(\alpha), \Im(\beta) \ll X^{1-\varepsilon}$ with an error term $O(X^{1/2+\varepsilon})$. The absence of any constraint on the imaginary parts of $\alpha$ and $\beta$ is a noteworthy advantage of Theorem \ref{Theorem for all characters}. \newline

   As a special case of Theorem \ref{Theorem for all characters}, one may consider the case $\beta \rightarrow \infty$. This amounts to a replacement of the expression $L(\tfrac 12+\beta, \chi^{(c_K\varpi)})$  by $1$ in \eqref{Integral for all characters} below. One then proceed as in Section \ref{sec 3} by  observing that as the parameter $\beta$ is not involved anymore, one may discard any dependence on $\beta$ in this process. This leads an expression equivalent to the limit case $\beta \rightarrow \infty$ on both sides of \eqref{Asymptotic for ratios of all characters} with the factor involving with $\beta$ in the error term being discarded.  We thus obtain the following result concerning the first moment of quadratic Hecke $L$-functions.  We also note here that it may be possible to make the following theorem unconditional in a manner similar to the proof of \cite[Theorem 1.2]{G&Zhao14}.
\begin{theorem}
\label{Thmfirstmoment}
With the notation as above and assuming the truth of GRH, we have, for $-1/4<\Re(\alpha)<1/2$,
\begin{align} \label{Asymptotic for first moment}
\begin{split}			
 \sum_{\substack{\varpi \odd}}\Lambda_{K}(\varpi) L(\tfrac 12+\alpha, \chi^{(c_K\varpi)})w \bfrac {N(\varpi)}X  
=&   X\M w(1)\zeta^{(2)}_K(1+2\alpha)+X^{1-\alpha}\M w(1-\alpha)B^{-\alpha}_K\Gamma_K(\tfrac 12+\alpha)\zeta^{(2)}_K(1-2\alpha) \\
& \hspace*{1.2cm} +O\lz|1+\alpha|^{3-d_K/2+\varepsilon}X^{\max (1/2, 3/4-\Re(\alpha))+\varepsilon}\pz. 
\end{split}
\end{align}
\end{theorem}

As the error term in \eqref{Asymptotic for first moment} is uniform in $\alpha$, taking $\alpha \rightarrow 0^+$ readily leads to the following asymptotic formula for the first moment of central values of the family of quadratic Hecke characters of prime-related moduli.  As our main focus is the error term, we omit the explicit expression of $Q$ appearing in the main term here. 
\begin{corollary}
\label{Thmfirstmomentatcentral}
		With the notation as above and assume the truth of GRH, we have,
\begin{align*}
%%\label{Asymptotic for first moment at central}
\begin{split}			
	& 	\sum_{\substack{\varpi \odd}}\Lambda_{K}(\varpi)L(\tfrac 12, \chi^{(c_K\varpi)})w \bfrac {N(\varpi)}X   = XQ(\log X)+O\lz  X^{3/4+\varepsilon}\pz.
\end{split}
\end{align*}
  where $Q$ is a linear polynomial whose coefficients depend only on the absolute constants $\M w(1)$ and $\M w'(1)$.
\end{corollary}

   Similarly, we take $\alpha \rightarrow \infty$ on both sides of \eqref{Asymptotic for ratios of all characters}, noting that we do not encounter the pole at $s=1-\alpha$ in this case. Hence, after dropping the second main term in \eqref{Asymptotic for ratios of all characters} and the factor involving with $\alpha$ in the error term, we obtain the following result concerning the negative first moment of quadratic Hecke $L$-functions.
\begin{theorem}
\label{Thmnegfirstmoment}
		With the notation as above and assuming the truth of GRH. We have for $\Re(\beta)>0$ and any $\varepsilon>0$,
\begin{align*}
%%\label{Asymptotic for neg first moment}
\begin{split}			
& \sum_{\substack{\varpi \odd}}\frac {\Lambda_{K}(\varpi)}{L(\tfrac 12+\beta, \chi^{(c_K\varpi)})}w \bfrac {N(\varpi)}X  =  X\M w(1)
 +O\lz |1+\beta|^{\varepsilon} X^{\max(1/2, 1-\Re(\beta))+\varepsilon}\pz.
\end{split}
\end{align*}
\end{theorem}

   In section \ref{sec: logder}, we shall also differentiate with respect to $\alpha$ in \eqref{Asymptotic for ratios of all characters} and set $\alpha=\beta=r$ to arrive at an asymptotic formula for the smoothed first moment of $L'(\frac{1}{2}+r,\chi^{(c_K\varpi)})/L(\frac{1}{2}+r,\chi^{(c_K\varpi)})$.
\begin{theorem}
\label{Theorem for log derivatives}
	With the notation as above and assuming the truth of GRH, we have for $1/4< \Re(r)<1/2$,
\begin{align}
\label{Sum of L'/L with removed 2-factors}
\begin{split}
  \sum_{\substack{\varpi}}\frac {\Lambda_{K}(\varpi)L'(\tfrac 12+r, \chi^{(c_K\varpi)})}{L(\tfrac 12+r, \chi^{(c_K\varpi)})}w \bfrac {N(\varpi)}X = X \M w(1) & \frac{(\zeta^{(j)}_K(1+2r))'}{\zeta^{(j)}_K(1+2r)}-X^{1-r} \M w(1-r)\Gamma_K(\tfrac 12+r)\frac {\zeta^{(2)}_K(1-2r)}{r_KB^{r}_K} \\
&+O(|1+r|^{1/2+\varepsilon}X^{1-\Re(r)+\varepsilon}).
\end{split}
\end{align}
\end{theorem}

Theorem \ref{Theorem for log derivatives} allows one to compute the one-level density of low-lying zeros of the corresponding families of quadratic Hecke $L$-functions, following the approach in the proof of \cite[Corollary 1.5]{Cech1} using \cite[Theorem 1.4]{Cech1}.  To this end, let $h(x)$ be an even Schwartz function such that its Fourier transform is supported in the interval $[-a,a]$ for some $a>0$.  We define the one-level density of the family of $L$-functions considered in this paper by
\begin{align*}
%%\label{Ddef}
\begin{split}
		D_K(X;h)=\frac{1}{F_K(X)}\sum_{\substack{ \varpi \odd}}\Lambda_K(\varpi)w \bfrac {N(\varpi)}X\sum_{\gamma_{\varpi, n}}h\bfrac{\gamma_{\varpi, n}\log X}{2\pi},
\end{split}
\end{align*}
   where $\gamma_{\varpi, n}$ runs over the imaginary parts of the non-trivial zeros of $L(s,\chi^{(c_K\varpi)})$ and
\begin{align}
\label{Size of family}
\begin{split}
			F_K(X)&=\sum_{\substack{\varpi \odd}}\Lambda_K(\varpi) w \bfrac {N(\varpi)}X.
\end{split}
\end{align}
	
Theorem \ref{Theorem for log derivatives} now enables us to computes $D_K(X;h)$ asymptotically.
\begin{theorem}
\label{Theorem one-level density}
	With the notation as above and assuming the truth of GRH, for any function $w(t)$ that is non-negative and compactly supported on the set of positive real numbers, we have
\begin{align}
\label{Onelevel}
\begin{split}
		D_K(X;h) 
=&\frac {2\widehat h(1)}{F_K(X) \mathcal{L}}  \sum_{\substack{\varpi \odd}}\Lambda_K(\varpi) w \bfrac {N(\varpi)}X\log N(\varpi) \\
& +\frac{2}{F_K(X) \mathcal{L}} \int\limits_{-\infty}^{\infty}h(u)\lz X \M w(1) \frac{(\zeta^{(j)}_K(1+ \frac{4\pi iu}{\mathcal{L}}))'}{\zeta^{(j)}_K(1+\frac{4\pi iu}{\mathcal{L}} )}-X^{1-2\pi iu/\mathcal{L}} \M w(1-\tfrac{4\pi iu}{\mathcal{L}})\Gamma_K(\tfrac 12+\tfrac {2\pi iu}{\mathcal{L}})\frac {\zeta^{(2)}_K(1-\frac {4\pi iu}{\mathcal{L}})}{r_KB_K^{2\pi iu/\mathcal{L}}} \pz \dif u \\
&+\frac {d_K}{2 \mathcal{L}}\int\limits^{\infty}_{-\infty}h\lz u \pz \lz \frac {\Gamma'}{\Gamma}(\tfrac {d_K}4+\tfrac {d_K \pi iu}{\mathcal{L}}) +\frac {\Gamma'}{\Gamma}(\tfrac {d_K}4-\tfrac {d_K \pi iu}{\mathcal{L}}) \pz \dif u+\frac {2\widehat h(1)}{\mathcal{L}} \log B_K +O(\frac{X^{(1+a)/2+\varepsilon}}{F_K(X)}),
\end{split}
\end{align}
where $\mathcal{L} = \log X$.  In particular, if $a<1$, then
\begin{align}
\label{Onelevelasym}
\begin{split}
		D_K(X;h)=& \int\limits^{\infty}_{-\infty}h(x) \dif x  +\frac{2\widehat h(1)}{\mathcal{L}} \lz \frac {1}{\widehat w(1)}\int\limits^{\infty}_0 w(u)\log u \dif u+d_K\frac {\Gamma'}{\Gamma}(\frac {d_K}4) +\log B_K \pz +O\left( \frac 1{\mathcal{L}^2} \right).
\end{split}
\end{align}
\end{theorem}

  Our proof of Theorem \ref{Theorem for all characters} uses the powerful method of multiple Dirichlet series, which has been introduced to study the first moment of the quadratic Dirichlet $L$-functions attached to fundamental discriminants by D. Goldfeld and J. Hoffstein 
in \cite{DoHo}, who improved the error term obtained in the asymptotic formula in an earlier work of M. Jutila  \cite{Jutila}. A systematic development on the method of multiple Dirichlet series can be found in \cite{DGH}.  Additionally, our work here also drew much inspiration from M. \v{C}ech \cite{Cech1}.  While typical usage of the multiple Dirichlet series often seeks to obtain analytical continuation of the underlying series to the whole complex space,  our approach differs from this treatment in the way that we only aim to obtain meromorphic continuation of the series involved to a region large enough that suffices to complete our task. This {\it modus operandi} has the advantage of not requiring adjustments to the multiple Dirichlet series with extra correction factors, often necessary in the usual approach.

\section{Preliminaries}
\label{sec 2}

\subsection{Primary Elements}
\label{sect: Kronecker}

   For an imaginary quadratic number field $K=\mq(\sqrt{-d})$ with $d \in \mathcal  S$, we cite the following facts about $K$ from \cite[Section 3.8]{iwakow}. \newline

  The ring of integers $\mathcal{O}_K$ is a free $\mz$ module (see \cite[Section 3.8]{iwakow}) such that $\mathcal{O}_K=\mz+\omega_K \mz$, where
\begin{align*}
   \omega_K & =\begin{cases}
\displaystyle     \frac {1+\sqrt{d}}{2} \qquad & d \equiv 1 \pmod 4, \\ \\
     \sqrt{d} \qquad & d \equiv 2, 3 \pmod 4.
    \end{cases}
\end{align*}
   The discriminants are given by
\begin{align*}
   D_K & =\begin{cases}
     d \qquad & \text{if $d \equiv 1 \pmod 4$}, \\
     4d \qquad & \text{if $d \equiv 2, 3 \pmod 4$}.
    \end{cases}
\end{align*}
  Note further that $D_K=1$ for $K=\mq$. \newline

   It is well-known that when $K$ is an imaginary quadratic number field of class number one or $K=\mq$, every ideal in $\mathcal O_K$ is principal so that one may fix a unique generator for each non-zero ideal.  We shall only work with ideals that are co-prime to $(2)$ in this paper and now determine, for each such ideal, a unique generator which we call the primary element.  We remark here that the definition for these primary elements is based on a result of F. Lemmermeyer \cite[Theorem 12.17]{Lemmermeyer05}.  Additionally, we shall not work directly with the concrete form of these primary elements since what matters to us is that Lemmas \ref{Quadrecgenprim} and \ref{lem: primquadGausssum} below are applicable to them.  \newline

   For $K=\mq$, we define an odd rational integer $n$ to be primary if $n>0$. For $K=\mq(\sqrt{d})$ with $d \in \mathcal S$, such primary elements are already given in \cite[Section 2.2]{G&Zhao16}, so we simply quote the definition from there. For $K=\mq(i)$, we define for any $n \in \mathcal O_K$ with $(n, 2)=1$ to be primary if and only if $n \equiv 1 \pmod {(1+i)^3}$. 
   For other values of $d$, we write $d=4k+1$ or $d=-2$ and we define
\begin{align*}
%%\label{G2square}
\begin{split}
   G=
\begin{cases}
 \{1, \omega^2_K=k+\omega_K, (1+\omega_K)^2=k+1-\omega_K \}, & 2 \nmid k, \\
  \{1 \}, & 2|k, \\
 \{1, -1+2\omega_K \} , & d=-2.
\end{cases}
\end{split}
\end{align*}
 
  We then define the primary elements to be the ones that are 
\begin{align*}
%%\label{primary}
\begin{split}
   \equiv
\begin{cases}
 G \times <1+2\omega_K> \pmod 4, & d \neq -2 \; \mbox{or} \; -3,  \\
  G \times \{1, -(1+\omega_K)\}  \pmod 4, & d = -2, \\
 <1+2\omega_K> \pmod 4, & d=-3.
\end{cases}
\end{split}
\end{align*}

%%----------------------------------------------------------------------------
\subsection{Quadratic characters and quadratic Gauss sums}
\label{sec2.4}
%%----------------------------------------------------------------------------
  Let $K$ be any imaginary quadratic number field or $K=\mq$. For an prime $\varpi \in \mathcal{O}_{K}$ that is co-prime to $2$, the quadratic symbol $\leg {\cdot}{\varpi}$ is defined for $a \in \mathcal{O}_{K}$, $(a, \varpi)=1$ by $\leg{a}{\varpi} \equiv a^{(N(\varpi)-1)/2} \pmod{\varpi}$, with $\leg{a}{\varpi} \in \{\pm 1 \}$.  If $\varpi | a$, we define $\leg{a}{\varpi} =0$. The quadratic symbol is extended to $\leg {\cdot}{n}$ for any $n$ co-prime to $2$ multiplicatively. We further define $\leg {\cdot}{c}=1$ for $c \in U_K$. \newline

 We note the following quadratic reciprocity law concerning primary elements from \cite[Proposition 1.11]{Lemmermeyer} for the case $K=\mq$ and from 
 \cite[Lemma 2.4, (2.7), (2.12)]{G&Zhao18}. 
\begin{lemma}
\label{Quadrecgenprim}
   Let $K=\mq(\sqrt{d})$ with $d \in \mathcal{S}$ or $K=\mq$. For any co-prime primary elements $n,m$ with $(nm, 2)=1$, we have
\begin{align*}
%%\label{quadreciKprim}
    \leg {n}{m}\leg{m}{n}=(-1)^{(N(n)-1)/2\cdot (N(m)-1)/2}.
\end{align*}
\end{lemma}

  Recall the definition of $c_K$ in Section \ref{sec 1} and we define $\chi^{(c_Kn)}$ to be the quadratic symbol $\leg {c_Kn}{\cdot}$ for any primary $n \in \mathcal O_K$.  Since $\chi^{(c_Kn)}$ is trivial on the units, we shall henceforth regard it as a Hecke character of trivial infinite type. We then define
the associated Gauss sum $g_K(\chi^{(c_Kn)})$ to be
\begin{align*}
%%\label{g2}
 g_K(\chi^{(c_Kn)}) := \sum_{x \shortmod{q}} \chi^{(c_Kn)}(x) \widetilde{e}_K\leg{x}{q}, \quad \mbox{where} \quad   \widetilde{e}_K(z) :=e\left( \text{\it Tr}\big(\frac {z}{\sqrt{D_K}}\big) \right).
\end{align*}
Here $e(z) = \exp (2 \pi i z)$ for any complex number $z$ and $\text{\it Tr}(n)$ denotes the trace of any $n \in K$. \newline

  Our next lemma evaluates $g_K(\chi^{(c_Kn)})$ and is quoted from \cite[Theorem 9.17]{MVa1} for the case $K=\mq$ and from 
 \cite[Lemma 2.6]{G&Zhao16} for the case $K$ that is an imaginary quadratic number field of class number one. 
\begin{lemma}
\label{lem: primquadGausssum}
  Let $K=\mq(\sqrt{d})$ with $d \in \mathcal S$ or $K=\mq$. For any odd, square-free $c \in \mathcal{O}_K$, we have
\begin{align*}
%%\label{primquadGausssum}
 g_K(\chi^{(c_Kc)})=N(c_Kc )^{1/2}.
\end{align*}
\end{lemma}

\subsection{Hecke $L$-functions}
For $K \neq \ratq$ and any primitive Hecke character $\chi$ of trivial infinite type, a well-known result of E. Hecke gives that $L(s, \chi)$ has an
analytic continuation to the whole complex plane and satisfies the functional equation (see \cite[Theorem 3.8]{iwakow})
\begin{align}
\label{fneqn}
  \Lambda(s, \chi) = W(\chi)\Lambda(1-s, \overline \chi), \; \mbox{where} \;  \big |W(\chi)\big | =1, 
\end{align}
and
\begin{align}
\label{Lambda}
  \Lambda(s, \chi) = (|D_K|N(q))^{s/2}(2\pi)^{-s}\Gamma(s)L(s, \chi).
\end{align}

  When $\chi=\chi^{(c_K\varpi)}$ for a primary prime $\varpi$, this functional equations can be written as 
\begin{align}
\label{fneqnL}
\begin{split}
  L(s, \chi^{(c_K\varpi)})=& g_{K}(c_K\varpi)N(c_K\varpi)^{-1/2}(|D_K|N(c_K\varpi ))^{1/2-s}(2\pi)^{2s-1}\frac {\Gamma(1-s)}{\Gamma (s)}L(1-s, \chi^{(c_K\varpi)}) \\
  =& (B_KN(\varpi))^{1/2-s}\Gamma_K(s)L(1-s, \chi^{(c_K\varpi)}),
\end{split}
\end{align}
 where $B_K$ and $\Gamma_K$ are defined in \eqref{GammaKBKdef} and the last equality above follows from Lemma \ref{lem: primquadGausssum}. \newline
 
  Similarly, the functional equation for Dirichlet $L$-functions given in \cite[\S 9]{Da} together with Lemma  \ref{lem: primquadGausssum} for the value of quadratic Gauss sum imply that $L(s, \chi^{(c_K\varpi)})$ also equals the last expression in \eqref{fneqnL} in the case $K=\mq$. \newline
  
 To estimate ratios of Gamma functions, we note that Stirling's formula (see \cite[(5.113)]{iwakow}) yields for constants $c_0, d_0 \in \mr$,
\begin{align}
\label{Stirlingratio}
\begin{split}
  \frac {\Gamma(c_0(1-s)+ d_0)}{\Gamma (c_0s+ d_0)} \ll (1+|s|)^{c_0(1-2\Re (s))}. 
\end{split}
\end{align}

  For any imaginary quadratic number field of class number one or $\mq$ and any Hecke character $\chi$ of trivial infinite type, the associated $L$-function $L(s, \chi)$ has an Euler product for $\Re(s)$ large enough given by
\begin{align*}
%%\label{LEuler}
 L(s, \chi)=\prod_{(\varpi)}\Big(1-\frac {\chi(\varpi)}{N(\varpi)^s}\Big)^{-1},
\end{align*}
  where $(n)$ denotes the ideal generated by any $n \in \mathcal O_K$. Logarithmically differentiating both sides above then implies that for $\Re(s)$ large enough,
\begin{align*}
%%\label{Llogder}
 -\frac {L'(s, \chi)}{L(s, \chi)}=\sum_{(n)}\frac {\Lambda_{K}(n)\chi(n)}{N(n)^s}.
\end{align*}
  Also, upon writing $\mu_K$ for the M\"obius function on $\mathcal O_K$, we have for $\Re(s)$ large enough,
\begin{align*}
%%\label{Llogder}
 \frac {1}{L(s, \chi)}=\sum_{(n)}\frac {\mu_{K}(n)\chi(n)}{N(n)^s}.
\end{align*}
  Similar to \cite[(2.9), (2.12)--(2.14)]{G&Zhao17}, we have that under GRH, for any primary $n$,
\begin{align} \label{Lderboundgen}
\begin{split}
  (s-1)\cdot -\frac {L'(s, \chi^{(c_Km)})}{L(s, \chi^{(c_Km)})}  \ll  |s-1|\big((N(m)+2)(|1+s|)\big)^{\varepsilon}, & \quad \Re(s) \geq 1/2+\varepsilon , \\
  (s-1)L(s, \chi^{(c_Km)}) \ll  |s-1|\big((N(m)+2)(|1+s|)\big)^{\varepsilon}, & \quad \Re(s) \geq 1/2 , \\
  L(s,  \chi^{(c_Km)})^{-1} \ll  |sN(m)|^{\varepsilon}, & \quad \Re(s) \geq 1/2+\varepsilon.
\end{split}
\end{align}

\subsection{Some results on multivariable complex functions} We shall require some results from multivariable complex analysis. We begin with the definition of a tube domain.
\begin{defin}
		An open set $T\subset\mc^n$ is a tube if there is an open set $U\subset\mr^n$ such that $T=\{z\in\mc^n:\ \Re(z)\in U\}.$
\end{defin}
	
   For a set $U\subset\mr^n$, we define $T(U)=U+i\mr^n\subset \mc^n$.  We shall make use of the following Bochner's Tube Theorem \cite{Boc}.
\begin{theorem}
\label{Bochner}
		Let $U\subset\mr^n$ be a connected open set and $f(z)$ a function holomorphic on $T(U)$. Then $f(z)$ has a holomorphic continuation to the convex hull of $T(U)$.
\end{theorem}

 We denote the convex hull of an open set $T\subset\mc^n$ by $\widehat T$.  Our next result is \cite[Proposition C.5]{Cech1} on the modulus of holomorphic continuations of multivariable complex functions.  The upshot is that analytic extensions inherit the same bound of the function from which they emanate.
\begin{prop} \label{Extending inequalities}
Assume that $T\subset \mc^n$ is a tube domain, $g,h:T\rightarrow \mc$ are holomorphic functions, and let $\tilde g,\tilde h$ be their holomorphic continuations to $\widehat T$. If  $|g(z)|\leq |h(z)|$ for all $z\in T$ and $h(z)$ is nonzero in $T$, then also $|\tilde g(z)|\leq |\tilde h(z)|$ for all $z\in \widehat T$.
\end{prop}

\section{Proof of Theorem \ref{Theorem for all characters}}
\label{sec 3}
\subsection{Setup}
   
   As the proofs are similar, we consider only the case for $K$ being an imaginary quadratic number field of class number one here. The Mellin inversion yields that
\begin{equation} \label{Integral for all characters}
	\sum_{\substack{\varpi \odd}}\frac {\Lambda_{K}(\varpi) L(\tfrac 12+\alpha, \chi^{(c_K\varpi)})}{L(\tfrac 12+\beta, \chi^{(c_K\varpi)})}w \bfrac {N(\varpi)}X=\frac1{2\pi i}\int\limits_{(c)}A_K\lz s,\tfrac12+\alpha, \tfrac12+\beta \pz X^s\widehat w(s) \dif s,
\end{equation}
  where, for $\Re(s), \Re(w)$ and $\Re(z)$ large enough,
\begin{align} \label{Aswzexp}
\begin{split}
A_K(s,w, z)=& \sum_{\substack{\varpi \odd}}\frac{\Lambda_{K}(\varpi) L(w,  \chi^{(c_K\varpi)})}{ N(\varpi)^sL(z,  \chi^{(c_K\varpi)})}=\sum_{\substack{ n \odd}}\frac{\Lambda_{K}(n) L(w, \chi^{(c_Kn)})}{N(n)^sL(z,  \chi^{(c_Kn)})}-\sum_{i \geq 2}\sum_{\substack{\varpi \odd}}\frac{\Lambda_{K}(\varpi^i) L(w, \chi^{(c_K\varpi^i)})}{N(\varpi)^{is}L(z,  \chi^{(c_K\varpi^i)})}. 
\end{split}
\end{align}
 Here we recall that the Mellin transform $\widehat{w}$ of $w$ is defined in \eqref{Mellin} and we also recall we use $\varpi$ to denote a prime in $K$ throughout this paper.

\subsection{Analytical properties of $A_K(s,w,z)$}
  
Note that 
\begin{align} 
\label{Aswremainder1}
\begin{split}
\sum_{i \geq 2} & \sum_{\substack{ \varpi \odd}} \frac{\Lambda_{K}(\varpi^i) L(w, \chi^{(c_K\varpi^i)})}{N(\varpi)^{is}L(z,  \chi^{(c_K\varpi)})} \\
=& \sum_{i \geq 1}\sum_{\substack{ \varpi \odd}}\frac{\Lambda_{K}(\varpi) L(w, \chi^{(c_K\varpi^{2i})})}{N(\varpi)^{2is}L(z,  \chi^{(c_K\varpi^{2i})})}+\sum_{i \geq 1}\sum_{\substack{\varpi \odd}}\frac{\Lambda_{K}(\varpi) L(w, \chi^{(c_K\varpi^{2i+1})})}{N(\varpi)^{(2i+1)s}L(z,  \chi^{(c_K\varpi^{2i+1})})} \\
=& \sum_{i \geq 1}\sum_{\substack{\varpi \odd}}\frac{\Lambda_{K}(\varpi) L(w, \chi^{(c_K)})(1-\frac {\chi^{(c_K)}(\varpi)}{N(\varpi)^w})}{N(\varpi)^{2is}L(z,  \chi^{(c_K)})(1-\frac {\chi^{(c_K)}(\varpi)}{N(\varpi)^z})}+\sum_{i \geq 1}\sum_{\substack{\varpi \odd}}\frac{\Lambda_{K}(\varpi) L(w, \chi^{(c_K\varpi)})}{N(\varpi)^{(2i+1)s}L(z,  \chi^{(c_K\varpi)})} =: \Sigma_1 + \Sigma_2,
\end{split}
\end{align}
say.  It is readily seen that $\Sigma_1$ converges in the region
\begin{equation*}
%%\label{key}
		S_{1,1}=\{(s,w,z): \ \Re(s)> \tfrac 12, \Re(z) > \tfrac 12 \}.
\end{equation*} 
 
  In what follows,  we shall define similar regions $S_{i, j}$ and $S_j$ and adopt the convention that for any real number $\delta$,
\begin{align*}
%%\label{Aswboundwlarge}
\begin{split}
 S_{i, j,\delta} :=& \{ (s,w,z)+\delta (1,1,1) : (s,w,z) \in S_{i,j}, \ \tfrac 14 \leq \Re(w) \leq 1 \}, \quad \mbox{and} \\
 S_{j,\delta} :=& \{ (s,w,z)+\delta (1,1,1) : (s,w,z) \in S_j, \ \tfrac 14 \leq \Re(w) \leq 1 \}.
\end{split}
\end{align*}

From \eqref{Lderboundgen}, the functional equation \eqref{fneqnL} by noting that it is still valid for $L(w, \chi^{(c_K)})$ upon replacing $\varpi$ by $1$ throughout and the estimation \eqref{Stirlingratio} that under GRH, we have in the region $S_{1,1, \varepsilon}$, 
\begin{align*}
%%\label{Aswboundwlarge1}
\begin{split}
\Sigma_1 \ll & |1+w|^{\frac 12+\varepsilon}|z|^{\varepsilon}.
\end{split}
\end{align*}   
   
   Further, $\Sigma_2$ converges under GRH for $\Re(s) > 1/3, \Re(w) \geq 1/2$ by \eqref{Lderboundgen}. Applying the functional equation \eqref{fneqnL}  for $L(w, \chi^{(c_K\varpi)})$ renders that $\Sigma_2$ is also convergent, under GRH, for $\Re(w)<1/2, \Re(3s+w)>3/2$. It then follows from Theorem \ref{Bochner} that $\Sigma_2$ converges in
\begin{equation*}
%%\label{key}
		S_{1,2}=\{(s,w,z): \ \Re(s)> \tfrac 13, \ \Re(3s+w)> \tfrac 32 , \ \Re(z) > \tfrac12 \}.
\end{equation*} 
 
     Similar to our discussions above, we have under GRH that in the region $S_{1,2, \varepsilon}$,
\begin{align*}
%%\label{Aswboundwlarge2}
\begin{split}
\Sigma_2 \ll  |1+w|^{\frac 12+\varepsilon}|z|^{\varepsilon}.
\end{split}
\end{align*}        
     
     We then conclude that the left-hand side expression in \eqref{Aswremainder1} is convergent in the region
\begin{equation*}
%%\label{key}
		S_{1,3} :=S_{1,1} \cap S_{1,2}=\{(s,w, z): \ \Re(s)> \tfrac 12, \ \Re(3s+w)> \tfrac 32 , \ \Re(z) > \tfrac12 \}.
\end{equation*}     
    Also, in the region $S_{1, 3, \varepsilon}$, under GRH,
\begin{align}
\label{Aswboundwlarge3}
\begin{split}
 \Big |\sum_{i \geq 2}\sum_{\substack{ \varpi \odd}}\frac{\Lambda_{K}(\varpi^i) L(w, \chi^{(c_K\varpi^i)})}{N(\varpi)^{is}L(z, \chi^{(c_K\varpi^i)})}\Big | \ll & |1+w|^{\frac 12+\varepsilon}|z|^{\varepsilon}.
\end{split}
\end{align}         

Next note that we have, using $\chi_{mk}$ to denote the character $\leg{\cdot}{mk}$,
\begin{align} 
\label{Aswmain}
\begin{split}
 \sum_{\substack{n \odd}}\frac{\Lambda_{K}(n) L(w, \chi^{(c_Kn)})}{N(n)^sL(z, \chi^{(c_Kn)})}  =& \sum_{\substack{n \odd}} \sum_{\substack{m, k \odd}}\frac{\mu_K(k)\Lambda_{K}(n) \chi_{mk}(c_Kn)}{N(m)^wN(n)^sN(k)^z} \\
 = & \sum_{\substack{m, k \odd }}\frac{\mu_K(k)\chi_{mk}(c_K)}{N(m)^wN(k)^z}\sum_{\substack{n \odd}}\frac{\Lambda_{K}(n)\chi_{mk}(n)}{N(n)^s} \\
=& \sum_{\substack{m, k \odd }}\frac{\mu_K(k)\chi_{mk}(c_K)}{N(m)^wN(k)^z}\sum_{\substack{ n \odd}}\frac{\Lambda_{K}(n)\chi^{(
(-1)^{(N(mk)-1)/2}mk)}(n)}{N(n)^s},
\end{split}
\end{align}
  where the last equality above follows from the quadratic reciprocity law. \newline

 We further recast the last expression above as
\begin{align} 
\label{Aswmainmcongcond}
\begin{split}
 &  \sum_{\substack{ m, k \odd \\ N(mk) \equiv 1 \pmod 4}} \frac{\mu_K(k)\chi_{mk}(c_K)}{N(m)^wN(k)^z}\sum_{\substack{ n \odd}}\frac{\Lambda_{K}(n)\chi^{(4mk)}(n)}{N(n)^s}+\sum_{\substack{m, k \odd \\ N(mk) \equiv -1 \pmod 4}}\frac{\chi_m(c_K)}{N(m)^w}\sum_{ n \odd}\frac{\Lambda_{K}(n)\chi^{(-4mk)}(n)}{N(n)^s} \\
=& \sum_{\substack{m, k \odd \\ N(mk) \equiv 1 \pmod 4}}\frac{\mu_K(k)\chi_{mk}(c_K)}{N(m)^wN(k)^z}\Big( -\frac {L'(s, \chi^{(4mk)})}{L(s, \chi^{(4mk)})} \Big) +\sum_{\substack{m, k \odd \\ N(mk) \equiv -1 \pmod 4}}\frac{\mu_K(k)\chi_{mk}(c_K)}{N(m)^wN(k)^z} \Big( -\frac {L'(s, \chi^{(-4mk)})}{L(s, \chi^{(-4mk)})} \Big). 
\end{split}
\end{align}

Now from \eqref{Lderboundgen} that under GRH, other than a simple pole at $s=1$ in the case of $mk$ being a perfect square, the last expression in \eqref{Aswmainmcongcond} is convergent in the region when $\Re(s)>1/2, \Re(w)>1, \Re(z)>1$. On the other hand,
\begin{align} \label{Aswmainalter}
 \sum_{\substack{ n \odd}}\frac{\Lambda_{K}(n) L(w, \chi^{(c_Kn)})}{N(n)^s L(z, \chi^{(c_Kn)})}=\sum_{n \odd}\frac{\Lambda_{K}(n)\chi^{(4)}(n) L(w, \chi^{(c_Kn)})}{N(n)^s L(z, \chi^{(c_Kn)})}.
\end{align}
  As the right-hand side expression of \eqref{Aswmainalter} is convergent when $\Re(s)>1$, $\Re(w)\geq 1/2$, $\Re(z)>1/2$ under GRH by \eqref{Lderboundgen}, it follows from \eqref{Aswmain}--\eqref{Aswmainalter} and Theorem \ref{Bochner} that the expression in \eqref{Aswmain} is convergent in the region
\begin{equation*}
%%\label{key}
		S_{1} :=\{(s,w, z): \ \Re(s)> \tfrac 12, \ \Re(w) \geq \tfrac 12, \ \Re(z) > \tfrac 12, \ \Re(s+w)>\tfrac 32, \ \Re(s+z)>\tfrac 32  \}.
\end{equation*} 
Furthermore, from \eqref{Lderboundgen} and Proposition \ref{Extending inequalities}, we deduce that, under GRH, in the region $S_{1,\varepsilon}$, 
\begin{align}
\label{Aswboundwlarge5}
\begin{split}
 \Big|(s-1)\sum_{\substack{n \odd}}\frac{\Lambda_{K}(n) L(w, \chi^{(c_Kn)})}{N(n)^sL(z, \chi^{(c_Kn)})}\Big | \ll & |1+s|^{1+\varepsilon}|wz|^{\varepsilon}.
\end{split}
\end{align}       

   Observe that $S_{1,3}$ contains $S_{1}$, so that by our discussions above and \eqref{Aswzexp} that $A_K(s,w,z)$ is analytic in the region $S_{1}$. It also follows from \eqref{Aswboundwlarge3} and \eqref{Aswboundwlarge5} that under GRH,  in the region $S_{1, \varepsilon}$,
\begin{align}
\label{Aswboundwlarge6}
\begin{split}
 \Big|(s-1)A_K(s,w,z) \Big | \ll & |1+s|^{1+\varepsilon}|1+w|^{\varepsilon}|z|^{\varepsilon}.
\end{split}
\end{align}    
   
  Lastly, we apply the functional equation \eqref{fneqnL} for $L(w, \chi^{(c_K\varpi)})$ to deduce from \eqref{Aswzexp} that
\begin{align} 
\label{Aswfuneqn}
\begin{split}
A_K(s,w,z)=& \sum_{\substack{\varpi \odd}}\frac{\Lambda_{K}(\varpi) L(w,  \chi^{(c_K\varpi)})}{N(\varpi)^sL(z,  \chi^{(c_K\varpi)})} 
= B_K^{1/2-w}\Gamma_K(w) \sum_{\substack{ \varpi \odd}}\frac{\Lambda_{K}(\varpi) L(1-w,  \chi^{(c_K\varpi)})}{N(\varpi)^{s+w-1/2}L(z,  \chi^{(c_K\varpi)})} \\
=& B_K^{1/2-w}\Gamma_K(w) A_K(s+w-\tfrac{1}{2}, 1-w, z). 
\end{split}
\end{align}  

  It follows from our discussions above that $A_{K}(s+w-1/2, 1-w, z)$ and hence $A_{K}(s,w, z)$ is convergent in the region
\begin{equation*}
%%\label{key}
		S_{2} :=\{(s,w, z): \ \Re(s+w)>1, \ \Re(1-w) \geq \tfrac 12, \ \Re(z)> \tfrac 12, \ \Re(s)> 1, \ \Re(s+w+z)> 2 \}.
\end{equation*}

  We now deduce from \eqref{Aswboundwlarge6}, \eqref{Aswfuneqn} and \eqref{Stirlingratio} that under GRH, in the region $S_{2, \varepsilon}$, we have for any $\varepsilon>0$, 
\begin{align}
\label{Aswboundwlarge7}
\begin{split}
 \Big|(s+w-3/2)A_K(s,w,z) \Big | \ll & |1+w|^{\frac 12+\varepsilon}|1+s+w|^{1+\varepsilon}|2-w|^{\frac 12+\varepsilon}|z|^{\varepsilon} \\
\ll & |1+s|^{1+\varepsilon}|1+w|^{2+\varepsilon}|z|^{\varepsilon}.
\end{split}
\end{align}    

  Note that the union of $S_1$ and $S_2$ is connected and the points $(1/2, 1)$, $(1, 0)$ are on the boundary of this union. The convex hull of $S_1 \cup S_2$ is
\begin{equation*}
%%\label{key}
		S_3 :=\{(s,w,z): \ \Re(s)>\tfrac 12, \ \Re(z)>\tfrac 12, \ \Re(s+w+z)> 2, \ \Re(2s+w+z)>3, \ \Re(2s+w)>2,\ \Re(s+w)>1,\ \Re(s+z)> \tfrac{3}{2}\}.
\end{equation*} 
 
    To see this, note that the two planes: $\Re(s+w)=3/2$, $\Re(s+z)= 3/2$ intersect when $\Re(z)=1/2$ at the line given by the intersection of $\Re(s+w)=3/2$ and $\Re(s)= 1$. Note also that the plane $\Re(s+w+z)= 2$ can be written as $\Re(s+w)= 2-\Re(z)$ and $2-\Re(z)$ equals $3/2$ when $\Re(z)=1/2$. Note also that $2-\Re(z)=1$ when $\Re(z)=1$. It follows that the convex hull of $S_1$ and $S_2$ contains the plane determined by the two lines: $\ \Re(s+w)=3/2$, $\Re(s+z)=3/2$ on $S_1$ and $\ \Re(s+w+z)=2,\ \Re(s)=1$ on $S_2$, intersecting with the two planes: $\Re(z)=1/2$ and $\Re(z)=1$. Note that the two lines both intersect the plane $\Re(z)=1/2$ at the point $(1, 1/2, 1/2)$. The first line intersects the plane $\Re(z)=1$ at $(1/2, 1, 1)$, the second line intersects the plane $\Re(z)=1$ at $(1, 0, 1)$. These three points then determine a plane whose equation is given by: $\Re(2s+w+z)=3$. Moreover, when $\Re(z)>1$, the convex hull continues with the plane that is determined by the lines: $\Re(s)=1/2$, $\Re(s+w)=3/2$ in $S_1$ and $\Re(s)=1, \ \Re(s+w)=1$ in $S_2$. We then see easily that the plane has the equation: $\Re(2s+w)=2$. \newline

  It follows from Theorem \ref{Bochner} that $A_K(s,w,z)$ has meromorphic continuation in the region $S_3$. Furthermore, we deduce from Proposition~\ref{Extending inequalities} the bound, inherited from \eqref{Aswboundwlarge6} and \eqref{Aswboundwlarge7}, that in the region $S_{3,\varepsilon}$, we have
\begin{align}
\label{Aswboundslarge2}
\begin{split}
 |(s+w-3/2)(s-1)A_{K}(s,w,z)| \ll & |1+s|^{2+\varepsilon}|1+w|^{2+\varepsilon}|z|^{\varepsilon}.
\end{split}
\end{align}

\subsection{Residues}
\label{sec:resA}

  We see from \eqref{Aswmainmcongcond} that $A_{K}(s,w,z)$ has a simple pole at $s=1$ arising from the terms with $mk$ being perfect squares.  As the residue of $-\zeta'_K(s)/\zeta_K(s)$ at $s=1$ is $1$, we deduce  that the residue at $s=1$ of the left-hand side expression in \eqref{Aswmainmcongcond} equals
\begin{align}
\label{Ress=1}
\begin{split}
  \sum_{\substack{m,k \odd \\ mk= \square }}\frac{\mu_K(k)}{N(m)^wN(k)^z} 
 =&  \prod_{(\varpi, 2)=1}\Big(1+\Big(1-\frac 1{N(\varpi)^{z-w}} \Big)\frac 1{N(\varpi)^{2w}}\Big(1-\frac 1{N(\varpi)^{2w}} \Big)^{-1}\Big ) 
= \frac {\zeta^{(2)}_K(2w)}{\zeta^{(2)}_K(w+z)},
\end{split}
\end{align}
where $mk = \square$ means $mk$ is a perfect square.  Similarly, we see from \eqref{Aswfuneqn} that $A_{K}(s,w,z)$ has simple poles at $s+w=3/2$ and \eqref{Ress=1} yields 
\begin{align} 
\label{Ressw32}
\begin{split}
 \res_{s=3/2-w}A_{K}(s,w,z)=B_K^{1/2-w}\Gamma_K(w)\frac {\zeta^{(2)}_K(2-2w)}{\zeta^{(2)}_K(1-w+z)}. 
\end{split}
\end{align}  

\subsection{Completion of the proof}

 We shift the line of integration in \eqref{Integral for all characters} to $\Re(s)=E(\alpha,\beta)+\varepsilon$.  The integral on the new line can be absorbed into the $O$-term in \eqref{Asymptotic for ratios of all characters} upon using \eqref{Aswboundslarge2} and the observation that repeated integration by parts gives that, for any integer $E \geq 0$,
\begin{align*}
%%\label{whatbound}
 \widehat w(s)  \ll  (1+|s|)^{-E}.
\end{align*}
   We also encounter two simple poles at $s=1$ and $s=1-\alpha$ in the move with the corresponding residues given in \eqref{Ress=1} and \eqref{Ressw32}.  Direct computations now lead to the main terms given in \eqref{Asymptotic for ratios of all characters}, completing the proof of Theorem \ref{Theorem for all characters}.

\section{Proof of Theorem \ref{Theorem for log derivatives}}
\label{sec: logder}

	 Again we consider only the case for $K$ being an imaginary quadratic number field of class number one here. We first recast the expression in \eqref{Asymptotic for ratios of all characters} as
\begin{align}
\label{Asymptoticshortversion}
\begin{split}	
& \sum_{\substack{\varpi \odd}}\frac {\Lambda_{K}(\varpi) L(\tfrac 12+\alpha, \chi^{(c_K\varpi)})}{L(\tfrac 12+\beta, \chi^{(c_K\varpi)})}w \bfrac {N(\varpi)}X =   XM_1(\alpha,\beta)+X^{1-\alpha}M_2(\alpha,\beta)+R(X,\alpha,\beta),
\end{split}
\end{align}
   where $R(X,\alpha,\beta)$ is the $O$-term in \eqref{Asymptotic for ratios of all characters},
\begin{align*}
%%\label{M1M2}
	M_1(\alpha,\beta) := \M w(1)\frac {\zeta^{(2)}_K(1+2\alpha)}{\zeta^{(2)}_K(1+\alpha+\beta)}  \quad \mbox{and} \quad
 M_2(\alpha,\beta) := \M w(1-\alpha)B^{-\alpha}_K\Gamma_K(\tfrac 12+\alpha)\frac {\zeta^{(2)}_K(1-2\alpha)}{\zeta^{(2)}_K(1-\alpha+\beta)} .
\end{align*}
	Note that the expression on the left-hand side of \eqref{Asymptoticshortversion} and both $M_i(\alpha,\beta), i=1,2$ are analytic functions of $\alpha,\beta$.  Hence so is $R(X,\alpha,\beta)$. \newline
	
	We now differentiate the above terms with respect to $\alpha$ for a fixed $\beta=r$ with $1/4 <\Re(r) <1/2$, and then set $\alpha=\beta=r$.  We get
\begin{align}
\label{M1der}
\begin{split}
		\frac{\dif}{\dif \alpha} XM_1(\alpha,\beta)\Big\vert_{\alpha=\beta=r}
			&=X \M w(1)\frac{(\zeta^{(j)}_K(1+2r))'}{\zeta^{(j)}_K(1+2r)}.	
\end{split}
\end{align}

Similarly, due to the presence of the factor $(\zeta^{(2)}_K(1-\alpha+\beta))^{-1}$, only one term survives in the derivative of $X^ {1-\alpha}M_2(\alpha,\beta)$ with respect to $\alpha$ and we have
\begin{align}
\label{M2der}
\begin{split}
		\frac{\dif}{\dif \alpha}X^ {1-\alpha}M_2(\alpha,\beta)\Big\vert_{\alpha=\beta=r}	&=-X^{1-r} \M w(1-r)B_K^{-r}\Gamma_K(\tfrac 12+r)\frac {\zeta^{(2)}_K(1-2r)}{r_K},
\end{split}
\end{align}	
where $r_K$ denotes the residue of $\zeta_K(s)$ at $s=1$. \newline

Next, as $R(X,\alpha,\beta)$ is analytic in $\alpha$, Cauchy's integral formula renders
\begin{align*}
%%\label{Rder}
\begin{split}
		\frac{\dif }{\dif \alpha}R(X,\alpha,\beta)=\frac{1}{2\pi i}\int\limits_{C_\alpha}\frac{R(X,z,\beta)}{(z-\alpha)^2} \dif z,
\end{split}
\end{align*}
   where $C_\alpha$ is a circle centered at $\alpha$ of radius $\rho$ with $\varepsilon/2<\rho<\varepsilon$. As inspection of the proof of Theorem \ref{Theorem for all characters} given in Section \ref{sec 3} implies that we have $|R(X,z,\beta)| \ll |1+z|^{1/2+\varepsilon}|1+\beta|^{\varepsilon}$ when $\Re(z) \geq 1/4$. Consequently, for $\Re(\alpha) \geq 1/4+\varepsilon$, we have
\begin{align*}
%%\label{Rderest}
\begin{split}
\lab\frac{\dif}{\dif \alpha}R(X,\alpha,\beta)\rab\ll \frac1{\rho}\cdot\max_{z\in C_\alpha} |R(X,z,\beta)|\ll |1+\alpha|^{1/2+\varepsilon}|1+\beta|^{\varepsilon}X^{E(\alpha,\beta)+\varepsilon}.
\end{split}
\end{align*}

    We now set $\alpha=\beta=r$ and choose $\rho<\varepsilon<r-1/4$ to deduce from \eqref{Nab} and the above that
\begin{align}
\label{Rder1}
\begin{split}
		\lab\frac{\dif }{\dif \alpha}R(X,\alpha,\beta)\rab_{\alpha=\beta=r} \ll |1+r|^{1/2+\varepsilon}X^{1-\Re(r)+\varepsilon}.
\end{split}
\end{align}
	
   We derive \eqref{Sum of L'/L with removed 2-factors} from \eqref{Asymptoticshortversion}--\eqref{Rder1}, getting Theorem \ref{Theorem for log derivatives}.
	
\section{Proof of Theorem \ref{Theorem one-level density}}
\label{Section one-level density}
	
As before, we consider only the case of $K$ being an imaginary quadratic number field of class number one here.  We apply the residue theorem to arrive at, recalling that $\mathcal{L} = \log X$,
\begin{align}
\label{sumoverzeros}
\begin{split}
		\sum_{\gamma_{\varpi,n}}h\bfrac{\gamma_{\varpi,n}\mathcal{L}}{2\pi}=\frac{1}{2\pi i}\lz \ \int\limits_{(2)}-\int\limits_{(-1)} \ \pz h\lz\frac{\mathcal{L}}{2\pi i}(s-\tfrac{1}{2})\pz\frac {\Lambda'}{\Lambda}(s, \chi^{(c_K\varpi)}) \dif s.
\end{split}
\end{align}

Now the functional equation  \eqref{fneqn} gives that
\begin{align}
\label{Lambdarel}
\begin{split}
		\frac {\Lambda'}{\Lambda}(s, \chi^{(c_K\varpi)})=-\frac {\Lambda'}{\Lambda}(1-s,  \chi^{(c_K\varpi)}).
\end{split}
\end{align}

   It follows from \eqref{Lambda}, \eqref{sumoverzeros} and \eqref{Lambdarel} that
\begin{align*}
%%\label{sumoverzeros1}
\begin{split}
		\sum_{\gamma_{\varpi,n}}h\bfrac{\gamma_{\varpi,n}\mathcal{L}}{2\pi}=&\frac{2}{2\pi i}\int\limits_{(2)}h\lz\frac{\mathcal{L}}{2\pi i}(s-\tfrac{1}{2})\pz \lz \frac {\Lambda'}{\Lambda}(s, \chi^{(c_K\varpi)})\pz \dif s \\
=& \frac{1}{2\pi i}\int\limits_{(2)}h\lz\frac{\mathcal{L}}{2\pi i}(s-\tfrac{1}{2}) \pz \lz 2\frac {L'}{L}(s, \chi^{(c_K\varpi)})+2\frac {\Gamma'}{\Gamma}(s)+\log (B_KN(\varpi)) \pz \dif s.
\end{split}
\end{align*}

The above lead to
\begin{align}
\label{Dsimplified}
\begin{split}
	F_K(X)D_K(X;h) =& \frac{1}{2\pi i}\int\limits_{(2)}h\lz\frac{\mathcal{L}}{2\pi i}(s-1/2)\pz\sum_{\substack{\varpi \odd}}\Lambda_K(\varpi) w \bfrac {N(\varpi)}X \\
& \hspace{1in} \times \lz  2\frac {L'}{L}(s, \chi^{(c_K\varpi)})+2\frac {\Gamma'}{\Gamma}(s)+\log (B_KN(\varpi)) \pz \dif s.
\end{split}
\end{align}

From \cite[Lemma 10.1]{Cech1}, if the Fourier transform of $h$ is supported in the interval $[-a,a]$, then
\begin{equation}
\label{Lemma size of h}
		h\bfrac{s\mathcal{L}}{2\pi i}\ll \frac{X^{a\cdot\Re(s)}}{|s|^2(\mathcal{L})^2}.
\end{equation}

  We now evaluate
\begin{equation*}
	I:= \frac{2}{2\pi i}\int\limits_{(2)}h\lz\frac{\mathcal{L}}{2\pi i}(s-\tfrac{1}{2})\pz\sum_{\substack{\varpi \odd}}\Lambda_K(\varpi) w \bfrac {N(\varpi)}X 
\frac {L'}{L}(s, \chi^{(c_K\varpi)}) \dif s
\end{equation*}
  by setting $s=1/2+r+it$ with $0<r<1/2$ and applying \eqref{Sum of L'/L with removed 2-factors} to see that $I$ equals
\begin{align}
\label{Integral after ratios}
\begin{split}
\frac {1}{\pi} \int\limits_{-\infty}^{\infty} h\lz\frac{\mathcal{L}}{2\pi i}(r+it)\pz & \lz X \M w(1)\frac{(\zeta^{(j)}_K(1+2r+2it))'}{\zeta^{(j)}_K(1+2r+2it)}-X^{1-r-it} \M w(1-r-it)\Gamma_K(\tfrac 12+r+it)\frac {\zeta^{(2)}_K(1-2r-2it)}{r_KB_K^{r+it}}\pz \dif t \\
& +\frac{1}{2\pi i}\int\limits_{(1/2+r)}h\lz\frac{\mathcal{L}}{2\pi i}(s-\tfrac{1}{2})\pz O(|1+t|^{1/2+\varepsilon}X^{1-r+\varepsilon}) \dif t.
\end{split}
\end{align}
  We bound the last integral above by setting $r=1/2-\varepsilon$ and using \eqref{Lemma size of h} to see that it is
\begin{align*}
%%\label{Error in one-level density}
  \ll & X^{1-r+\varepsilon}\int\limits_{-\infty}^{\infty}h\lz\frac{\mathcal{L}}{2\pi i}(r+it)\pz|1+t|^{1/2+\varepsilon}|\dif t| \ll X^{1-r+ar+\varepsilon}, 
\end{align*}
 which gives the error term in \eqref{Onelevel}. \newline

 We also evaluate the first integral in \eqref{Integral after ratios} by shifting the line of integration to $r=0$.  The integral in question is
\begin{align}
\label{Main term in one-level density}
\begin{split}
& \frac {1}{\pi} \int\limits_{-\infty}^{\infty}h\lz \frac{\mathcal{L}}{2\pi}t\pz \lz  \frac{(\zeta^{(j)}_K(1+2it))'}{\zeta^{(j)}_K(1+2it)}-X^{1-it} \M w(1-it)\Gamma_K(\tfrac 12+it)\frac {\zeta^{(2)}_K(1-2it)}{r_KB_K^{it}}\pz \dif t\\
		=& \frac {2 }{ \mathcal{L}} \int\limits_{-\infty}^{\infty}h(u)\lz X \M w(1) \frac{(\zeta^{(j)}_K(1+ \frac {4\pi iu}{\mathcal{L}}))'}{\zeta^{(j)}_K(1+\frac {4\pi iu}{\mathcal{L}})}-X^{1-\frac {2\pi iu}{\mathcal{L}}} \M w(1-\tfrac {2\pi iu}{\mathcal{L}})\Gamma_K(\tfrac 12+\tfrac {2\pi iu}{\mathcal{L}})\frac {\zeta^{(2)}_K(1-\frac {4\pi iu}{\mathcal{L}})}{r_KB_K^{2\pi iu/\mathcal{L}}} \pz \dif u.
\end{split}
\end{align}

   We next compute
\begin{align*}
%%\label{Dsimplified1}
\begin{split}
	\frac{1}{2\pi i}\int\limits_{(2)}h\lz\frac{\mathcal{L}}{2\pi i}(s-\tfrac{1}{2})\pz\sum_{\substack{\varpi \odd}}\Lambda_K(\varpi) w \bfrac {N(\varpi)}X\lz 2\frac {\Gamma'}{\Gamma}(s)+\log (B_KN(\varpi)) \pz \dif s
\end{split}
\end{align*}
  by shifting the line of integration to $\Re(s)=1/2$ to see that it equals
\begin{align}
\label{Dsimplified2}
\begin{split}
 &	\frac {2\widehat h(1)}{\mathcal{L}} \sum_{\substack{\varpi \odd}}\Lambda_K(\varpi) w \bfrac {N(\varpi)}X\log (B_KN(\varpi))+ \sum_{\substack{\varpi \odd}}\Lambda_K(\varpi) w \bfrac {N(\varpi)}X \frac{1}{\pi}\int\limits^{\infty}_{-\infty}h\lz\frac{\mathcal{L}}{2\pi}t \pz \frac {\Gamma'}{\Gamma}(\tfrac 12+it)\dif t.
\end{split}
\end{align}

 As $h$ is even, we see that
\begin{align}
\label{Dsimplified3}
\begin{split}
	\frac{1}{\pi}\int\limits^{\infty}_{-\infty}h\lz\frac{\mathcal{L}}{2\pi}t \pz \frac {\Gamma'}{\Gamma}(\tfrac 12+it) \dif t=& \frac{1}{2\pi}\int\limits^{\infty}_{-\infty}h\lz\frac{\mathcal{L}}{2\pi}t \pz \lz \frac {\Gamma'}{\Gamma}(\tfrac 12+it) +\frac {\Gamma'}{\Gamma}(\tfrac 12-it) \pz \dif t \\
=& \frac 1{\mathcal{L}}\int\limits^{\infty}_{-\infty}h\lz u \pz \lz \frac {\Gamma'}{\Gamma}(\tfrac 12+\tfrac {2 \pi iu}{\mathcal{L}}) +\frac {\Gamma'}{\Gamma}(\tfrac 12-\tfrac {2 \pi iu}{\mathcal{L}}) \pz \dif u.
\end{split}
\end{align}

Now \eqref{Onelevel} readily follows from \eqref{Size of family}, \eqref{Dsimplified} and \eqref{Main term in one-level density}--\eqref{Dsimplified3}. We further simplify the right-hand side of \eqref{Onelevel}. For any real $y \geq 0$,
  it follows from \cite[(5.49)]{iwakow} and \cite[Theorem 5.15]{iwakow} that
\begin{align*}
%%\label{sumstarvarpiest}
\begin{split}
			 \sum_{\substack{\varpi \odd \\ N(\varpi) \leq y}}\Lambda_K(\varpi) =y+O(y^{1/2+\varepsilon}).
\end{split}
\end{align*}

   The above estimation and partial summation now yield
\begin{align}
\label{sumstarvarpilogest}
\begin{split}
		& F_K(X)= \int^{\infty}_0  w \bfrac {u}X  \dif \Big( \sum_{\substack{\varpi \odd \\ N(\varpi) \leq u}}\Lambda_K(\varpi) \Big)= \widehat w(1)X+O(X^{1/2+\varepsilon}),  \\
  &\sum_{\substack{\varpi \odd}}\Lambda_K(\varpi) w \bfrac {N(\varpi)}X\log N(\varpi) = \widehat w(1)X\mathcal{L}+\left( \int\limits^{\infty}_0w(u)\log u \dif u \right) X+O(X^{1/2+\varepsilon}).
\end{split}
\end{align}

  Note that for any $a, b \in \mr$, we have the approximate formula (cf.~\cite[8.363.3]{GR})
\begin{align*}
%%\label{Dsimplified4}
\begin{split}
  \frac{\Gamma'}{\Gamma}(a + ib) + \frac{\Gamma'}{\Gamma}(a-ib) = 2\frac{\Gamma'}{\Gamma} (a) + O \left( \leg ba^2
\right).
\end{split}
\end{align*}

Thus \eqref{Dsimplified3} is
\begin{align}
\label{Dsimplified5}
\begin{split}
 \frac{4 \widehat h(1)}{\mathcal{L}}\frac{\Gamma'}{\Gamma}(\tfrac 12)  + O(\mathcal{L}^{-3} ).
\end{split}
\end{align}

   We readily deduce \eqref{Onelevelasym} from \eqref{Onelevel}, \eqref{sumstarvarpilogest} and \eqref{Dsimplified5}. This completes the proof of Theorem \ref{Theorem one-level density}.

\vspace*{.5cm}

\noindent{\bf Data Availability Statement .}  This manuscript has no associated data. \newline

\noindent{\bf Acknowledgments.}  P. G. is supported in part by NSFC grant 12471003 and L. Z. by the FRG PS71536 at the University of New South Wales.  The authors would like to thank the anonymous referee for his/her very careful inspection of the paper and many helpful suggestions.

\bibliography{biblio}
\bibliographystyle{amsxport}

\end{document}